\documentclass[12pt]{article}
\usepackage[english]{babel}
\usepackage[cp1251]{inputenc}
\begin{document}
\begin{center}

{\bf DECOMPOSITIONS FOR HYPERGEOMETRIC FUNCTION  \\
$ H_A \left( {\alpha ,\beta _1 ,\beta _2 ;\gamma _1 ,\gamma _2
;x,y,z} \right),$ $ H_B \left( {\alpha ,\beta _1 ,\beta _2 ;\gamma
_1 ,\gamma ,\gamma ;x,y,z} \right),$ $H_C \left( {\alpha ,\beta _1
,\beta _2 ;\gamma \,;x,y,z} \right)$ }\\[2mm]
\medskip
{\bf Anvar H. Hasanov, Rakhila B. Seilkhanova, Roza D. Seilova}\\
\medskip
S. Baishev University Aktobe, Department for Information Systems and Applied Mathematics, \\
030000 Aktobe, Zhubanov str. 302, Kazakhstan \\
E-Mails: anvarhasanov@yahoo.com, srahila@inbox.ru, roza$_{-}$seilova@mail.ru
\medskip
\end{center}
\medskip

{\bf Abstract.} With the help of some techniques based on certain
inverse pairs of symbolic operators, the authors investigated
several decomposition formulas associated with Srivastava's
Hypergeometric functions of three variables. Some operator
identities have been constructed in this matter. With the help
these operator forms 15 decompositions are found which are
expressed through product of Hypergeometric Gauss and Appell's
functions.



{\bf 1. Introduction.}
\medskip

A great interest to the theory of Hypergeometric functions of
several variables was due to Solutions of many applied problems of
differential equations of partial derivatives which are used with
the help of such Hypergeometric functions. For instance, the
energy absorbed with nonferromagnet conductor sphere included in
internal magnet field, are calculated with the help of such
functions [1]. The Hypergeometric functions of several variables
are used in quantized chemistry as well [2]. Especially, in gas
dynamics many problems lead to solution of degenerating second
order partial differential equations. For instance, the problem of
adiabatic flat-parallel without whirlwind gas flow and the flow
problem of supersonic current from vessel with flat walls, also a
lot of problems connecting with gas flow [3].

We note that the Riemann's functions and fundamental solutions of
the degenerating second order partial differential equations are
being defined by Hypergeometric functions of several variables
[4]. In investigation the boundary value problems for these
equations we need decompositions for Hypergeometric functions of
several variables which are expressed through simple
Hypergeometric Gauss and Appell's functions.

The operator method has been used in papers [5-7] and found
decompositions for Hypergeometric functions of two variables which
are expressed through Gauss Hypergeometric functions of one
variables.

Over four decades ago by H.M.Srivastava was introduced
Hypergeometric functions of three variables [8]
$$
H_A \left( {\alpha ,\beta _1 ,\beta _2 ;\gamma _1 ,\gamma _2
;x,y,z} \right) = \sum\limits_{m,n,p = 0}^\infty  {} \frac{{\left(
\alpha  \right)_{m + p} \left( {\beta _1 } \right)_{m + n} \left(
{\beta _2 } \right)_{n + p} }}{{\left( {\gamma _1 } \right)_m
\left( {\gamma _2 } \right)_{n + p} m!n!p!}}x^m y^n z^p ,\eqno
(1.1)
$$

$$
H_B \left( {\alpha ,\beta _1 ,\beta _2 ;\gamma _1 ,\gamma _2
,\gamma _3 ;x,y,z} \right) = \sum\limits_{m,n,p = 0}^\infty  {}
\frac{{\left( \alpha  \right)_{m + p} \left( {\beta _1 }
\right)_{m + n} \left( {\beta _2 } \right)_{n + p} }}{{\left(
{\gamma _1 } \right)_m \left( {\gamma _2 } \right)_n \left(
{\gamma _3 } \right)_p m!n!p!}}x^m y^n z^p,\eqno(1.2)
$$

$$ H_C \left( {\alpha ,\beta _1 ,\beta _2 ;\gamma \,;x,y,z} \right)
= \sum\limits_{m,n,p = 0}^\infty  {} \frac{{\left( \alpha
\right)_{m + p} \left( {\beta _1 } \right)_{m + n} \left( {\beta
_2 } \right)_{n + p} }}{{\left( \gamma  \right)_{m + n + p}
m!n!p!}}x^m y^n z^p.\eqno (1.3)
$$

In this work we construct decompositions for Hypergeometric
functions (1.1)-(1.2) with the help of Burchnall and Chaundy [5-7]
method. These decompositions consist of simple Appell's functions
and Gauss functions. Further, by means of obtained decompositions
we determine some definite integrals connecting with functions
(1.1)-(1.2).\\[2mm]

{\bf 2. Main operators.}
\medskip

Over six decades ago, Burchnall and Chaundy [5, 6] and Chaundy [7]
systematically presented a number of expansion and decomposition
formulas for double hypergeometric functions in series of simplest
Hypergeometric functions. Their method is based on  following
inverse pairs of symbolic operators:

$$
\nabla _{xy} \left( h \right): = \frac{{\Gamma \left( h
\right)\Gamma \left( {\delta _1  + \delta _2  + h}
\right)}}{{\Gamma \left( {\delta _1  + h} \right)\Gamma \left(
{\delta _2  + h} \right)}} = \sum\limits_{k = 0}^\infty
{\frac{{\left( { - \delta _1 } \right)_k \left( { - \delta _2 }
\right)_k }}{{\left( h \right)_k k!}}},\eqno (2.1)$$

$$\Delta _{xy} \left( h \right): = \frac{{\Gamma \left( {\delta _1
 + h} \right)\Gamma \left( {\delta _2  + h} \right)}}{{\Gamma
 \left( h \right)\Gamma \left( {\delta _1  + \delta _2  + h} \right)}}
 = \sum\limits_{k = 0}^\infty  {\frac{{\left( { - \delta _1 } \right)_k
 \left( { - \delta _2 } \right)_k }}{{\left( {1 - h - \delta _1  - \delta _2
 } \right)_k k!}}}$$

$$
= \sum\limits_{k = 0}^\infty  {\frac{{\left( { - 1}
  \right)^k \left( h \right)_{2k} \left( { - \delta _1 } \right)_k \left( {
   - \delta _2 } \right)_k }}{{\left( {h + k - 1} \right)_k \left( {\delta _1
    + h} \right)_k \left( {\delta _2  + h} \right)_k k!}}},\eqno   (2.2)
$$

$$\nabla _{xy} \left( h \right)\Delta _{xy} \left( l \right): = \frac{{\Gamma
 \left( h \right)\Gamma \left( {\delta _1  + \delta _2  + h} \right)}}{{\Gamma
 \left( {\delta _1  + h} \right)\Gamma \left( {\delta _2  + h} \right)}}\frac{{\Gamma
  \left( {\delta _1  + l} \right)\Gamma \left( {\delta _2  + l} \right)}}{{\Gamma
  \left( l \right)\Gamma \left( {\delta _1  + \delta _2  + l} \right)}}$$
$$= \sum\limits_{k = 0}^\infty  {\frac{{\left( {l - h} \right)_k \left( l
\right)_{2l} \left( { - \delta _1 } \right)_k \left( { - \delta _2
} \right)_k }}{{\left( {l + k - 1} \right)_k \left( {l + \delta _1
} \right)_k \left( {l + \delta _2 } \right)_k k!}}}$$
$$= \sum\limits_{k = 0}^\infty
 {\frac{{\left( {h - l} \right)_k \left( l \right)_{2l} \left( { - \delta _1 }
 \right)_k \left( { - \delta _2 } \right)_k }}{{\left( h \right)_k \left( {1 -
 l - \delta _1  - \delta _2 } \right)_k k!}}},\,\,\,\,(\delta _1 : = x\frac{\partial
 }{{\partial x}};\,\,\delta _2 : = y\frac{\partial }{{\partial y}}). \eqno  (2.3)$$

We now introduce here the following multivariable analogues of the
Burchnall- Chaundy symbolic operators $ \nabla _{xy} \left( h
\right)$ and $ \Delta _{xy} \left( h \right)$ defined by (2.1) -
(2.3), respectively (cf. [9, p. 240]; see [10, p. 113] for the
case when $ r = 3$ ):

$$ \tilde \nabla _{x_1 ;x_2  \cdot  \cdot  \cdot x_r } \left( h
\right): = \frac{{\Gamma \left( h \right)\Gamma \left( {\delta _1
+  \cdot  \cdot  \cdot  + \delta _r  + h} \right)}}{{\Gamma \left(
{\delta _1  + h} \right)\Gamma \left( {\delta _2  +  \cdot  \cdot
\cdot  + \delta _r  + h} \right)}} =$$ $$= \sum\limits_{k_2
,...,k_r = 0}^\infty  {} \frac{{\left( { - \delta _1 }
\right)_{k_2  + ... + k_r } \left( { - \delta _2 } \right)_{k_2 }
\cdot  \cdot  \cdot \left( { - \delta _r } \right)_{k_r }
}}{{\left( h \right)_{k_2  + ... + k_r } k_2 ! \cdot  \cdot  \cdot
k_r !}}, \eqno (2.4)$$

and

$$\tilde \Delta _{x_1 ;x_2  \cdot  \cdot  \cdot x_r } \left( h \right): =
\frac{{\Gamma \left( {\delta _1  + h} \right)\Gamma \left( {\delta _2  +
 \delta _3  +  \cdot  \cdot  \cdot  + \delta _r  + h} \right)}}{{\Gamma
 \left( h \right)\Gamma \left( {\delta _1  + \delta _2  +  \cdot  \cdot
 \cdot  + \delta _r  + h} \right)}}$$ $$= \sum\limits_{k_2 ,...,k_r  =
 0}^\infty  {} \frac{{\left( { - \delta _1 } \right)_{k_2  + ... + k_r }
 \left( { - \delta _2 } \right)_{k_2 }  \cdot  \cdot  \cdot \left( { -
 \delta _r } \right)_{k_n } }}{{\left( {1 - h - \delta _1  - \delta _2
 -  \cdot  \cdot  \cdot  - \delta _r } \right)_{k_2  + ... + k_r } k_2 !
 \cdot  \cdot  \cdot k_r !}}$$ $$ = \sum\limits_{k_2 ,...,k_r  = 0}^\infty
{} \frac{{\left( { - 1} \right)^{k_2  + ... + k_r } \left( \gamma  \right)_{2\left(
{k_2  + ... + k_r } \right)} \left( { - \delta _1 } \right)_{k_2  + ... + k_r }
\left( { - \delta _2 } \right)_{k_2 }  \cdot  \cdot  \cdot \left( { - \delta _n }
 \right)_{k_r } }}{{\left( {\gamma  + k_2  + ... + k_r  - 1} \right)_{k_2  + ... +
  k_r } \left( {\delta _1  + \gamma } \right)_{k_2  + ... + k_r } \left( {\delta _2
   + \delta _3  +  \cdot  \cdot  \cdot  + \delta _r  + \gamma } \right)_{k_2
   + ... + k_r } k_2 ! \cdot  \cdot  \cdot k_r !}} , \eqno  (2.5)$$

$$ (\delta _j : = x_j \frac{\partial }{{\partial x_j }};\,\,j =
1,...,r),$$ where we have applied such known multiple
hypergeometric summation as (cf. [11]; see also [12, p. 117])

$$F_D^{\left( r \right)} \left( {a;b_1 , \cdot  \cdot  \cdot ,b_r
;c;1, \cdot  \cdot  \cdot ,1} \right) = \frac{{\Gamma \left( c
\right)\Gamma \left( {c - a - b_1  -  \cdot  \cdot  \cdot  - b_r }
\right)}}{{\Gamma \left( {c - a} \right)\Gamma \left( {c - b_1  -
\cdot  \cdot  \cdot  - b_r } \right)}}, \eqno (2.6)$$
$$({\mathop{\rm Re}\nolimits} \,(c - a - b_1  -  \cdot  \cdot \cdot
- b_r ) > 0;\,c \notin Z_0 : = \left\{ {0, - 1, - 2, - 3, \cdot
\cdot \cdot } \right\}),$$ for the Lauricella function $
F_D^{\left( r \right)}$ in $ r$ variables, defined by (cf. [11];
see also [13, p.33, Equation 1.4 (4)])

$$F_D^{\left( r \right)} \left( {a;b_1 , \cdot  \cdot  \cdot ,b_r
;c;x_1 ,x_2 , \cdot  \cdot  \cdot ,x_r } \right): =$$ $$:=
\sum\limits_{m_1 ,m_2 , \cdot  \cdot  \cdot ,m_r  = 0}^\infty  {}
\frac{{\left( a \right)_{m_1  + m_2  +  \cdot  \cdot  \cdot  + m_r
} \left( {b_1 } \right)_{m_1 }  \cdot  \cdot  \cdot \left( {b_r }
\right)_{m_r } }}{{\left( c \right)_{m_1  + m_2  +  \cdot  \cdot
\cdot  + m_r } m_1 !m_2 ! \cdot  \cdot  \cdot m_r !}}x_1^{m_1 }
x_2^{m_2 }  \cdot  \cdot  \cdot x_r^{m_r },\eqno  (2.7)$$ $$\left(
{\max \left\{ {\left| {x_1 } \right|,...,\,\,\left| {x_r }
\right|} \right\} < 1} \right),$$ $(h)_k $ -  Pochhammer's symbol.\\[2mm]

{\bf 3. Main operator identities.}
\medskip

According to operators (2.1)-(2.5) for Hypergeometric functions
(1.1)-(1.3) we find the following operator identities

\medskip

$H_A \left( {\alpha ,\beta _1 ,\beta _2 ;\gamma _1 ,\gamma _2
;x,y,z} \right)$ $$= \nabla _{xz} \left( \alpha  \right)\nabla
_{xy} \left( {\beta _1 } \right)F\left( {\alpha ,\beta _1 ;\gamma
_1 ;x} \right)F_1 \left( {\beta _2 ;\beta _1 ,\alpha ;\gamma _2
;y,z} \right),\eqno  (3.1)$$

$H_B \left( {\alpha ,\beta _1 ,\beta _2 ;\gamma _1 ,\gamma ,\gamma
;x,y,z} \right) $ $$= \nabla _{xz} \left( \alpha  \right)\nabla
_{xy} \left( {\beta _1 } \right)\nabla _{yz} \left( \gamma
\right)F\left( {\alpha ,\beta _1 ;\gamma _1 ;x} \right)F_1 \left(
{\beta _2 ;\beta _1 ,\alpha ; \gamma ;y,z} \right), \eqno (3.2)$$

$H_C \left( {\alpha ,\beta _1 ,\beta _2 ;\gamma \,;x,y,z} \right)$
$$= \nabla _{xz} \left( \alpha  \right)\nabla _{xy} \left( {\beta _1 }
 \right)\tilde \Delta _x \left( \gamma  \right)F\left( {\alpha ,\beta _1
  ;\gamma \,;x} \right)F_1 \left( {\beta _2 ;\beta _1 ,\alpha ;\gamma
  \,;y,z} \right) ,\eqno  (3.3)$$

$ H_A \left( {\alpha ,\beta _1 ,\beta _2 ;\gamma _1 ,\gamma _2
;x,y,z} \right)$
$$= \nabla _{xz} \left( \alpha  \right)\nabla _{xy} \left( {\beta _1 } \right)
\Delta _{yz} \left( {\gamma _2 } \right)F\left( {\alpha ,\beta _1
;\gamma _1 ;x} \right)F_2 \left( {\beta _2 ;\beta _1 ,\alpha
;\gamma _2 ,\gamma _2 ;y,z} \right), \eqno (3.4)$$

$H_B \left( {\alpha ,\beta _1 ,\beta _2 ;\gamma _1 ,\gamma _2
,\gamma _3 ;x,y,z} \right) $
$$ = \nabla _{xz} \left( \alpha  \right)\nabla _{xy} \left( {\beta _1 }
 \right)F\left( {\alpha ,\beta _1 ;\gamma _1 ;x} \right)F_2 \left( {\beta
  _2 ;\beta _1 ,\alpha ;\gamma _2 ,\gamma _3 ;y,z} \right), \eqno (3.5)$$

$ H_C \left( {\alpha ,\beta _1 ,\beta _2 ;\gamma \,;x,y,z}
\right)$
$$  = \nabla _{xz} \left( \alpha  \right)\nabla _{xy} \left( {\beta _1 }
\right)\nabla _{yz} \left( \gamma  \right)\tilde \Delta _x \left( \gamma
 \right)F\left( {\alpha ,\beta _1 ;\gamma \,;x} \right)F_2 \left( {\beta
 _2 ;\beta _1 ,\alpha ;\gamma ,\gamma ;y,z} \right),\eqno   (3.6)$$

$ H_A \left( {\alpha ,\beta _1 ,\beta _2 ;\gamma _1 ,\gamma _2
;x,y,z} \right)$
$$  = \nabla _{xz} \left( \alpha  \right)\nabla _{xy} \left( {\beta _1 }
\right)\nabla _{yz} \left( {\beta _2 } \right)F\left( {\alpha
,\beta _1 ; \gamma _1 ;x} \right)F_3 \left( {\beta _1 ,\beta _2
,\beta _2 ,\alpha ; \gamma _2 ;y,z} \right),\eqno (3.7)$$

$ H_B \left( {\alpha ,\beta _1 ,\beta _2 ;\gamma _1 ,\gamma
,\gamma ;x,y,z} \right) $
$$  = \nabla _{xz} \left( \alpha  \right)\nabla _{xy} \left( {\beta _1 }
 \right)\nabla _{yz} \left( {\beta _2 } \right)\nabla _{yz} \left( \gamma
  \right)F\left( {\alpha ,\beta _1 ;\gamma _1 ;x} \right)F_3 \left(
  {\beta _2 ,\alpha ,\beta _1 ,\beta _2 ;\gamma ;y,z} \right) ,\eqno (3.8)$$

$ H_C \left( {\alpha ,\beta _1 ,\beta _2 ;\gamma \,;x,y,z} \right)
$ $$ = \nabla _{xz} \left( \alpha  \right)\nabla _{xy} \left(
{\beta _1 } \right)\nabla _{yz} \left( {\beta _2 } \right)\tilde
\Delta _x \left( \gamma \right)F\left( {\alpha ,\beta _1 ;\gamma
;x} \right)F_3 \left( {\beta _2 , \beta _2 ,\beta _1 ,\alpha
;\gamma ;y,z} \right),\eqno  (3.9)$$

$ H_A \left( {\alpha ,\alpha ,\beta _2 ;\gamma _1 ,\gamma _2
;x,y,z} \right) $ $$= \nabla _{xy} \left( \alpha  \right)\nabla _{xz}
\left( \alpha  \right)\Delta _{yz} \left( \alpha  \right)\Delta _{yz}
 \left( {\gamma _2 } \right)F\left( {\alpha ,\alpha ;\gamma _1 ;x}
 \right)F_4 \left( {\alpha ,\,\beta _2 ;\gamma _2 ,\gamma _2 ;y,z} \right) ,\eqno (3.10)$$

$ H_B \left( {\alpha ,\alpha ,\beta _2 ;\gamma _1 ,\gamma _2
,\gamma _3 ;x,y,z} \right) $ $$ = \nabla _{xz} \left( \alpha
\right)\nabla _{xy} \left( \alpha  \right)\Delta _{yz} \left(
\alpha  \right)F\left( {\alpha ,\alpha ; \gamma _1 ;x} \right)F_4
\left( {\alpha ,\beta _2 ,\gamma _2 ,\gamma _3 ;y,z} \right),\eqno
(3.11)$$

$ H_C \left( {\alpha ,\alpha ,\beta _2 ;\gamma \,;x,y,z} \right) $
$$  = \nabla _{xz} \left( \alpha  \right)\nabla _{xy} \left( \alpha
\right)\Delta _{yz} \left( \alpha  \right)\Delta _{yz} \left( \gamma
 \right)\tilde \Delta _x \left( \gamma  \right)F\left( {\alpha ,\alpha ;
 \gamma ;x} \right)F_4 \left( {\alpha ,\beta _2 ,\gamma ,\gamma ;y,z} \right), \eqno (3.12)$$

$ H_A \left( {\alpha ,\beta _1 ,\beta _2 ;\gamma _1 ,\gamma _2
;x,y,z} \right) $ $$= \nabla _{xz} \left( \alpha  \right)\nabla _{xy} \left(
 {\beta _1 } \right)\nabla _{yz} \left( {\beta _2 } \right)\Delta _{yz} \left(
  {\gamma _2 } \right)F\left( {\alpha ,\beta _1 ;\gamma _1 ;x} \right)F\left(
  {\beta _1 ,\beta _2 ;\gamma _2 ;y} \right)F\left( {\alpha ,\beta _2 ;\gamma _2
  ;z} \right) ,\eqno   (3.13)$$

$ H_B \left( {\alpha ,\beta _1 ,\beta _2 ;\gamma _1 ,\gamma _2
,\gamma _3 ;x,y,z} \right) $ $$ = \nabla _{xz} \left( \alpha  \right)\nabla _{xy}
 \left( {\beta _1 } \right)\nabla _{yz} \left( {\beta _2 } \right)F\left( {\alpha ,
 \beta _1 ;\gamma _1 ;x} \right)F\left( {\beta _1 ,\beta _2 ;\gamma _2 ;y} \right)F
 \left( {\alpha ,\beta _2 ;\gamma _3 ;z} \right) ,\eqno    (3.14)$$

$ H_C \left( {\alpha ,\beta _1 ,\beta _2 ;\gamma \,;x,y,z} \right)
$ $$= \nabla _{xz} \left( \alpha  \right)\nabla _{xy} \left(
{\beta _1 } \right) \nabla _{yz} \left( {\beta _2 } \right)\Delta
_{yz} \left( \gamma  \right)\tilde \Delta _x \left( \gamma
\right)F\left( {\alpha ,\beta _1 ;\gamma \,;x} \right)F\left(
{\beta _1 ,\beta _2 ;\gamma \,;y} \right)F\left( {\alpha ,\beta
_2; \gamma \,;z} \right) ,\eqno (3.15)$$

The operator identities for Hypergeometric functions $H_A \,,H_B
,H_C $ can be proved with the help of Mellin's transformation
[14]. For instance, we show that the operator identity (3.1) is
true. According to known E.W.Barnes definition [12, p.40] we have

$$F\left( {a;b;c;x} \right) = \frac{1}{{2\pi i}}\frac{{\Gamma
\left( c \right)}}{{\Gamma \left( a \right)\Gamma \left( b
\right)}}\int\limits_{ - i\infty }^{ + i\infty } {} \Gamma \left[
{\begin{array}{*{20}c}
   {a + s_1 ,b + s_1 , - s_1 }  \\
   {c + s_1 }  \\
\end{array}} \right]\left( { - x} \right)^{s_1 } ds_1 ,\eqno(3.16)$$
$c \ne 0, - 1, - 2,...$ and ${\mathop{\rm Re}\nolimits} \,s_1  <
0,\,{\mathop{\rm Re}\nolimits} \,\left( {a + s_1 } \right) >
0,\,{\mathop{\rm Re}\nolimits} \,\left( {b + s_1 } \right) > 0$
,\\[1mm]

$F_1 \left( {a;b_1 ,b_2 ;c;x,y} \right) =  - \frac{1}{{4\pi ^2 }}
\frac{{\Gamma \left( c \right)}}{{\Gamma \left( a \right)\Gamma \left( {b_1 }
 \right)\Gamma \left( {b_2 } \right)}}$ $$\cdot \int\limits_{ - i\infty }^{ +
 i\infty } {\int\limits_{ - i\infty }^{ + i\infty } {} } \Gamma
 \left[ {\begin{array}{*{20}c}
   {a + s_1  + s_2 ,b_1  + s_1 ,b_2  + s_2 , - s_1 , - s_2 }  \\
   {c + s_1  + s_2 }  \\
\end{array}} \right]\left( { - x} \right)^{s_1 } \left( { - y} \right)^{s_2 } ds_1
ds_2 , \eqno (3.17) $$ $c \ne 0, - 1, - 2,...$ and $${\mathop{\rm
Re}\nolimits} \,s_1  < 0,\,{\mathop{\rm Re}\nolimits} \,s_2  <
0,\,\,{\mathop{\rm Re}\nolimits} \,\left( {a + s_1  + s_2 }
\right) > 0,\,{\mathop{\rm Re}\nolimits} \,\left( {b_1  + s_1 }
\right) > 0,\,{\mathop{\rm Re}\nolimits} \,\left( {b_2  + s_2 }
\right) > 0,$$

By using equality $$ H_A \left( {\alpha ,\beta _1 ,\beta _2
;\gamma _1 ,\gamma _2 ;x,y,z} \right) = \sum\limits_{m = 0}^\infty
{} \frac{{\left( \alpha  \right)_m \left( {\beta _1 } \right)_m
}}{{\left( {\gamma _1 } \right)_m m!}}x^m F_1 \left( {\beta _2
;\beta _1  + m,\alpha + m;\gamma _2 ;y,z} \right),$$ we find
$$H_A \left( {\alpha ,\beta _1 ,\beta _2 ;\gamma _1 ,\gamma _2 ;x,y,z} \right)
=  - \frac{1}{{8\pi ^2 i}}\frac{{\Gamma \left( {\gamma _1 }
\right)\Gamma \left( {\gamma _2 } \right)}}{{\Gamma \left( \alpha
\right)\Gamma \left( {\beta _1 } \right)\Gamma \left( {\beta _2 }
\right)}}$$ $$ \cdot \int\limits_{ - i\infty }^{ + i\infty }
{\int\limits_{ - i\infty }^{ + i\infty } {\int\limits_{ - i\infty
}^{ + i\infty } {} } } \Gamma \left[ {\begin{array}{*{20}c}
   {\alpha  + s_1  + s_3 ,\beta _1  + s_1  + s_2 ,\beta _2  + s_2
    + s_3 , - s_1 , - s_2 , - s_3 }  \\
   {\gamma _1  + s_1 ,\gamma _2  + s_2  + s_3 }  \\
\end{array}} \right] $$
$$ \cdot \left( { - x} \right)^{s_1 } \left( { - y} \right)^{s_2 }
\left( { - z} \right)^{s_3 } ds_1 ds_2 ds_3  , \eqno (3.18)$$
$\gamma _1 \,,\,\gamma _2  \ne 0, - 1, - 2,...$ and ${\mathop{\rm
Re}\nolimits} \,s_1  < 0,$ ${\mathop{\rm Re}\nolimits} \,s_2  <
0,$ ${\mathop{\rm Re}\nolimits} \,s_2  < 0,$ ${\mathop{\rm
Re}\nolimits} \,\left( {\alpha  + s_1  + s_3 } \right) > 0,$
${\mathop{\rm Re}\nolimits} \,\left( {\beta _1  + s_1  + s_2 }
\right) > 0,\,{\kern 1pt} {\mathop{\rm Re}\nolimits} \,\left(
{\beta _2  + s_2  + s_3 } \right) > 0,$ where

$$\Gamma \left[ {\begin{array}{*{20}c}
   {a_1 ,a_2 ,...,a_m }  \\
   {b_1 ,b_2 ,...,b_n }  \\
\end{array}} \right] = \frac{{\Gamma \left( {a_1 } \right)\Gamma
\left( {a_2 } \right) \cdot  \cdot  \cdot \Gamma \left( {a_m }
\right)}}{{\Gamma \left( {b_1 } \right)\Gamma \left( {b_2 } \right)
 \cdot  \cdot  \cdot \Gamma \left( {b_n } \right)}}.$$

Therefore, from equalities (3.16)-(3.18) we can have the following

$$F\left( {a;b;c;x} \right) \leftrightarrow \frac{{\Gamma \left( c
\right)}}{{\Gamma \left( a \right)\Gamma \left( b \right)}}\Gamma
\left[ {\begin{array}{*{20}c}
   {a + s_1 ,b + s_1 , - s_1 }  \\
   {c + s_1 }  \\
\end{array}} \right],\eqno (3.19)$$ $c \ne 0, - 1, - 2,...$ and
${\mathop{\rm Re}\nolimits}\,s_1  < 0,\,{\mathop{\rm Re}\nolimits}
\,\left( {a + s_1 } \right) > 0,\,{\mathop{\rm Re}\nolimits}
\,\left( {b + s_1 } \right) > 0$ ,

$$F_1 \left( {a;b_1 ,b_2 ;c;x,y} \right) \leftrightarrow
\frac{{\Gamma \left( c \right)}}{{\Gamma \left( a \right)\Gamma
\left( {b_1 } \right)\Gamma \left( {b_2 } \right)}}\Gamma \left[
{\begin{array}{*{20}c}
   {a + s_1  + s_2 ,b_1  + s_1 ,b_2  + s_2 , - s_1 , - s_2 }  \\
   {c + s_1  + s_2 }  \\
\end{array}} \right],\eqno (3.20)$$ $c \ne 0, - 1, - 2,...$  and
${\mathop{\rm Re}\nolimits} \,s_1  < 0,$ ${\mathop{\rm
Re}\nolimits} \,s_2  < 0,$ ${\mathop{\rm Re}\nolimits} \,\left( {a
+ s_1  + s_2 } \right) > 0,$ ${\mathop{\rm Re}\nolimits} \,\left(
{b_1  + s_1 } \right) > 0,$ ${\mathop{\rm Re}\nolimits} \,\left(
{b_2  + s_2 } \right) > 0$ ,\\[1mm]

$H_A \left( {\alpha ,\beta _1 ,\beta _2 ;\gamma _1 ,\gamma _2
;x,y,z} \right) $ $$  \leftrightarrow \frac{{\Gamma \left( {\gamma
_1 } \right)\Gamma \left( {\gamma _2 } \right)}}{{\Gamma \left(
\alpha  \right)\Gamma \left( {\beta _1 } \right)\Gamma \left(
{\beta _2 } \right)}}\Gamma \left[ {\begin{array}{*{20}c}
   {\alpha  + s_1  + s_3 ,\beta _1  + s_1  + s_2 ,\beta _2  + s_2  + s_3 ,
   - s_1 , - s_2 , - s_3 }  \\
   {\gamma _1  + s_1 ,\gamma _2  + s_2  + s_3 }  \\
\end{array}} \right] ,\eqno(3.21)$$   $\gamma _1 \,,\,\gamma _2  \ne 0, - 1, - 2,...$
and ${\mathop{\rm Re}\nolimits} \,s_1  < 0,\,{\mathop{\rm
Re}\nolimits} \,s_2  < 0,\,\,{\mathop{\rm Re}\nolimits} \,s_2  <
0,\,\,$  ${\mathop{\rm Re}\nolimits} \,\left( {\alpha + s_1 + s_3
} \right) > 0,$ ${\mathop{\rm Re}\nolimits} \,\left( {\beta _1 +
s_1  + s_2 } \right) > 0,\,{\kern 1pt}$ $ {\mathop{\rm
Re}\nolimits} \,\left( {\beta _2  + s_2  + s_3 } \right) > 0$ .

At present we proof identity (3.1). Applying Mellin's
transformation in both parts of (3.1) and taking into account
equality $$\left( { - \delta _i } \right)_k f\left( x \right)
\leftrightarrow \left( { - s_i } \right)_k f^ *  \left( {s_i }
\right), \eqno (3.22)$$ $i = 1,2,3$ , moreover $ \mathop {\lim
}\limits_{x \to 0} \,x^{s - j - 1} f^{\left( j \right)} \left( x
\right) = 0,\,\,j = 0,1,...,k - 1,$ and from definition for (2.1)
we have\\[1mm]

$H_A \left( {\alpha ,\beta _1 ,\beta _2 ;\gamma _1 ,\gamma _2
;x,y,z} \right) $ $$ \leftrightarrow \frac{{\Gamma \left( {\gamma
_1 } \right)\Gamma \left( {\gamma _2 } \right)}}{{\Gamma ^2 \left(
\alpha  \right)\Gamma ^2 \left( {\beta _1 } \right)\Gamma \left(
{\beta _2 } \right)}}\Gamma \left[ {\begin{array}{*{20}c}
   {\alpha  + s_1 ,\alpha  + s_3 ,\beta _1  + s_1 ,\beta _1  + s_2
   ,\beta _2  + s_2  + s_3 , - s_1 , - s_2 , - s_3 ,}  \\
   {\gamma _1  + s_1 ,\gamma _2  + s_2  + s_3 }  \\
\end{array}} \right]$$ $$  \cdot \sum\limits_{i,j = 0}^\infty
{\frac{{\left( { - s_1 } \right)_i \left( { - s_1 } \right)_j
\left( { - s_2 } \right)_j \left( { - s_3 } \right)_i }}{{\left(
 \alpha  \right)_i \left( {\beta _1 } \right)_j i!j!}}}  \eqno (3.23)$$

By virtue of identity [15],

$$F\left( {a,b;c;1} \right) = \frac{{\Gamma \left( c \right)\Gamma
\left( {c - a - b} \right)}}{{\Gamma \left( {c - a} \right)\Gamma
\left( {c - b} \right)}},$$ $$({\mathop{\rm Re}\nolimits} \,(c - a
- b > 0,\,\,\,c \notin Z_0 : = \left\{ {0, - 1, - 2, - 3, \cdot
\cdot  \cdot } \right\}),$$ for expression (3.23) we have the
following $$ \sum\limits_{i,j = 0}^\infty  {\frac{{\left( { - s_1
} \right)_i \left( { - s_1 } \right)_j \left( { - s_2 } \right)_j
\left( { - s_3 } \right)_i }}{{\left( \alpha  \right)_i \left(
{\beta _1 } \right)_j i!j!}}}  = \frac{{\Gamma \left( \alpha
\right)\Gamma \left( {\alpha  + s_1  + s_3 } \right)}}{{\Gamma
\left( {\alpha  + s_1 } \right)\Gamma \left( {\alpha  + s_3 }
\right)}}\frac{{\Gamma \left( {\beta _1 } \right)\Gamma \left(
{\beta _1  + s_1  + s_2 } \right)}}{{\Gamma \left( {\beta _1  +
s_1 } \right)\Gamma \left( {\beta _1  + s_2 }
\right)}}.\eqno(3.24)$$

Substituting equalities (3.24) into expression (3.23) and after
simplest calculations we derive the operator identities (3.1).
Analogically, with the help of Mellin's transformation can be
proved the other operator identities.\\[2mm]

{\bf 4. The Decompositions for Hypergeometric functions $H_A,$
$H_B,$ $H_C$. }
\medskip

Applying operators and superposition of operators from identities
(3.1)-(3.15) we have the following
$$H_A \left( {\alpha ,\beta _1 ,\beta _2 ;\gamma _1 ,\gamma _2 ;x,y,z} \right) =
 \sum\limits_{i,j = 0}^\infty  {} \frac{{\left( \alpha  \right)_{i + j}
  \left( {\beta _1 } \right)_{i + j} \left( {\beta _2 } \right)_{i + j}
   }}{{\left( {\gamma _1 } \right)_{i + j} \left( {\gamma _2 } \right)_{i
    + j} i!j!}}x^{i + j} y^j z^i  \cdot$$ $$ \cdot F\left( {\alpha  + i + j,
    \beta _1  + i + j;\gamma _1  + i + j;x} \right)F_1 \left( {\beta _2  +
     i + j;\beta _1  + i + j,\alpha  + i;\gamma _2  + i + j;y,z} \right) ,\eqno (4.1)$$

$$ H_B \left( {\alpha ,\beta _1 ,\beta _2 ;\gamma _1 ,\gamma
,\gamma ;x,y,z} \right) = \sum\limits_{i,j,k = 0}^\infty  {}
\frac{{\left( \alpha  \right)_{i + j + k} \left( {\beta _1 }
\right)_{i + j + k} \left( {\beta _2 } \right)_{2i + j + k}
}}{{\left( {\gamma _1 } \right)_{j + k} \left( \gamma  \right)_i
\left( \gamma  \right)_{2i + j + k} i!j!k!}}x^{j + k} y^{i + j}
z^{i + k} $$ $$ F\left(
 {\alpha  + i + j + k,\beta _1  + i + j + k;\gamma _1  + j + k;x} \right) $$
$$ F_1 \left( {\beta _2  + 2i + j + k;\beta _1  + i + j,\alpha  + i + j + k;
\gamma  + 2i + j + k;y,z} \right), \eqno (4.2)$$

$$H_C \left( {\alpha ,\beta _1 ,\beta _2 ;\gamma \,;x,y,z} \right)
=$$ $$= \sum\limits_{i,j,k,l = 0}^\infty  {} \frac{{\left( { - 1}
\right)^{i
 + j} \left( \alpha  \right)_{i + 2j + k + l} \left( {\beta _1 } \right)_{2i
 + j + k + l} \left( {\beta _2 } \right)_{i + j + k + l} \left( \gamma  \right)_{2i
  + 2j} }}{{\left( {\gamma  + i + j - 1} \right)_{i + j} \left( \gamma  \right)_{2i
   + 2j + k + l} \left( \gamma  \right)_{2i + 2j + k + l} i!j!k!l!}}x^{i + j + k +
   l} y^{i + k} z^{j + l}  \cdot $$ $$\cdot F\left( {\alpha  + i + 2j + k + l,
   \beta _1  + 2i + j + k + l;\gamma  + 2i + 2j + k + l\,;x} \right) $$ $$ F_1 \left(
   {\beta _2  + i + j + k + l;\beta _1  + 2i + j + k,\alpha  + i + 2j + k + l;
   \gamma  + 2i + 2j + k + l\,;y,z} \right),\eqno  (4.3)$$

$$ H_A \left( {\alpha ,\beta _1 ,\beta _2 ;\gamma _1 ,\gamma _2
;x,y,z} \right) = $$ $$ = \sum\limits_{i,j,k = 0}^\infty  {\left(
{ - 1} \right)^k \frac{{\left( {\gamma _2 } \right)_{2k} \left(
\alpha  \right)_{i + j + k} \left( {\beta _1 } \right)_{i + j + k}
\left( {\beta _2 } \right)_{i + j + 2k} }}{{\left(
 {\gamma _2  + k - 1} \right)_k \left( {\gamma _1 } \right)_{i + j} \left(
 {\gamma _2 } \right)_{i + 2k} \left( {\gamma _2 } \right)_{j + 2k} i!j!k!}}}
  x^{i + j} y^{j + k} z^{i + k} $$

$$ F\left( {\alpha  + i + j + k,\beta _1
  + i + j + k;\gamma _1  + i + j;x} \right)$$ $$ F_2 \left( {\beta _2  + i + j
   + 2k;\beta _1  + j + k,\alpha  + i + j + k;\gamma _2  + j + 2k,\gamma _2  +
    i + 2k;y,z} \right) ,\eqno (4.4)$$

$$ H_B \left( {\alpha ,\beta _1 ,\beta _2 ;\gamma _1 ,\gamma _2
,\gamma _3 ;x,y,z}
 \right) = \sum\limits_{i,j = 0}^\infty  {} \frac{{\left( \alpha  \right)_{i + j}
 \left( {\beta _1 } \right)_{i + j} \left( {\beta _2 } \right)_{i + j} }}{{\left(
 {\gamma _1 } \right)_{i + j} \left( {\gamma _2 } \right)_j \left( {\gamma _3 }
 \right)_i i!j!}}x^{i + j} y^j z^i  \cdot $$ $$  \cdot F\left( {\alpha  + i + j,
 \beta _1  + i + j;\gamma _1  + i + j;x} \right)F_2 \left( {\beta _2  + i + j;
 \beta _1  + i + j,\alpha  + i;\gamma _2  + j,\gamma _3  + i;y,z} \right)\eqno (4.5)$$

$$ H_C \left( {\alpha ,\beta _1 ,\beta _2 ;\gamma \,;x,y,z} \right)
= $$ $$= \sum\limits_{i,j,k,l,r = 0}^\infty  {} \frac{{\left( { -
1} \right)^{i + j + r} \left( \alpha  \right)_{i + 2j + k + l + r}
\left( {\beta _1 } \right)_{2i + j + k + l} \left( {\beta _1 }
\right)_{2i + j + k + r} \left( {\beta _2 } \right)_{i + j + k + l
+ 2r} }}{{\left( {\gamma  + i + j - 1} \right)_{i + j} \left(
{\gamma  + 2i + 2j + k + l + r - 1} \right)_r \left( {\beta _1 }
\right)_{2i + j + k} }}$$ $$\frac{{\left( \gamma  \right)_{2i +
2j} }} {{\left( \gamma  \right)_{2i + 2j + k + l} \left( \gamma
\right)_{2i + 2j + k + l + 2r} i!j!k!l!r!}}x^{i + j + k + l} y^{i
+ k + r} z^{j + l + r}  \cdot  $$ $$\cdot F\left( {\alpha  + i +
2j + k + l, \beta _1  + 2i + j + k + l;\gamma  + 2i + 2j + k +
l\,;x} \right) \cdot  $$ $$F_2 (\beta _2  + i + j + k + l +
2r;\beta _1 + 2i + j + k + r,\alpha  + i + 2j + k + l + r; $$
$$\gamma  + 2i + 2j + k + l + 2r,\gamma  + 2i + 2j + k + l +
2r;y,z) ,\eqno (4.6)$$ $$H_A \left( {\alpha ,\beta _1 ,\beta _2 ;\gamma _1 ,
\gamma _2 ;x,y,z} \right) = \sum\limits_{i,j,k = 0}^\infty  {\frac{{\left(
\alpha  \right)_{i + j + k} \left( {\beta _1 } \right)_{i + j + k} \left(
{\beta _2 } \right)_{i + j} \left( {\beta _2 } \right)_{i + k} }}{{\left(
{\beta _2 } \right)_i \left( {\gamma _1 } \right)_{j + k} \left( {\gamma _2
 } \right)_{2i + j + k} i!j!k!}}} x^{j + k} y^{i + j} z^{i + k} $$
 $$ F\left( {\alpha  + i + j + k,\beta _1  + i + j + k;\gamma _1  + j + k;x}
 \right)$$ $$F_3 \left( {\beta _1  + i + j,\beta _2  + i + k,\beta _2  + i + j,
 \alpha  + i + j + k;\gamma _2  + 2i + j + k;y,z} \right) ,\eqno  (4.7)$$
$$ H_B \left( {\alpha ,\beta _1 ,\beta _2 ;\gamma _1 ,\gamma ,\gamma ;x,y,z} \right) = $$
$$= \sum\limits_{i,j,k,l = 0}^\infty  {\frac{{\left( \alpha  \right)_{j + 2k + l}
 \left( {\beta _1 } \right)_{i + j + k + l} \left( {\beta _2 } \right)_{2i + j +
 k} \left( {\beta _2 } \right)_{2i + j + l} }}{{\left( {\beta _2 } \right)_{2i +
  j} \left( \gamma  \right)_i \left( \gamma  \right)_{2i + 2j + k + l} \left(
  {\gamma _1 } \right)_{k + l} i!j!k!l!}}} x^{k + l} y^{i + j + k} z^{i + j + l}  \cdot $$
$$\cdot F\left( {\alpha  + j + 2k + l,\beta _1  + i + j + k + l;\gamma _1  + k +
l;x} \right) \cdot $$
$$\cdot F_3 ( \beta _2  + 2i + j + k,\alpha  + i + j + k + l,\beta _1
 + i + j + k, \beta _2  + 2i + j + l;$$ $$ \gamma  + 2i + 2j + k + l;y,z) ,\eqno (4.8)$$
$$H_C \left( {\alpha ,\beta _1 ,\beta _2 ;\gamma \,;x,y,z} \right) =  $$
$$ \sum\limits_{i,j,k,l,r = 0}^\infty  {} \frac{{\left( { - 1} \right)^{i + j}
 \left( \alpha  \right)_{i + 2j + k + l} \left( \alpha  \right)_{i + 2j + k + r}
  \left( {\beta _1 } \right)_{2i + j + k + l + r} \left( {\beta _2 } \right)_{i
   + j + k + l + r} \left( \gamma  \right)_{2i + 2j} }}{{\left( {\gamma  + i +
   j - 1} \right)_{i + j} \left( \alpha  \right)_{i + 2j + k} \left( \gamma
   \right)_{2i + 2j + k + l} \left( \gamma  \right)_{2i + 2j + k + l + 2r}
   i!j!k!l!r!}}\cdot$$ $$\cdot x^{i + j + k + l} y^{i + l + r} z^{j + k + r} \cdot
   F\left( {\alpha  + i + 2j + k + l,\beta _1  + 2i + j + k + l;\gamma
+ 2i + 2j + k + l;x} \right) $$ $$
 F_3 (\beta _2  + i + j + k + l + r,\beta _2  + i + j + k + l + r,\beta _1
   + 2i + j + k + l + r,\alpha  + i + 2j + k + r; $$
$$;\gamma  + 2i + 2j + k + l + 2r;y,z),\eqno      (4.9)$$

$$H_A \left( {\alpha ,\beta _1 ,\beta _2 ;\gamma _1 ,\gamma _2 ;x,y,z} \right) = $$
$$= \sum\limits_{i,j,k = 0}^\infty  {} \frac{{\left( \alpha  \right)_{i + j}
\left( \alpha  \right)_{i + k} \left( {\beta _1 } \right)_{i + j + k} \left(
 {\beta _2 } \right)_{i + j + k} \left( {\gamma _2  - \beta _2 } \right)_k }}
 {{\left( {\gamma _2  + i + j + k - 1} \right)_k \left( \alpha  \right)_i \left(
 {\gamma _1 } \right)_{i + j} \left( {\gamma _2 } \right)_{i + j + 2k} i!j!k!}}
 x^{i + j} y^{j + k} z^{i + k}  \cdot $$ $$ \cdot F\left( {\alpha  + i + j,\beta _1
  + i + j;\gamma _1  + i + j;x} \right) \cdot  $$ $$\cdot F\left( {\beta _2  + i
  + j + k,\beta _1  + i + j + k;\gamma _2  + i + j + 2k;y} \right) \cdot
  $$ $$\cdot F\left( {\beta _2  + i + j + k,\alpha  + i + k;\gamma _2  + i
  + j + 2k;z} \right) ,\eqno  (4.10)$$
$$ H_B \left( {\alpha ,\beta _1 ,\beta _2 ;\gamma _1 ,\gamma _2 ,\gamma _3 ;x,y,z}
 \right) = \sum\limits_{i,j,k = 0}^\infty  {} \frac{{\left( \alpha  \right)_{i + j
  + k} \left( {\beta _1 } \right)_{i + j + k} \left( {\beta _2 } \right)_{i + j}
  \left( {\beta _2 } \right)_{j + k} }}{{\left( {\beta _2 } \right)_j \left(
  {\gamma _1 } \right)_{i + k} \left( {\gamma _2 } \right)_{i + j} \left(
  {\gamma _3 } \right)_{j + k} i!j!k!}}x^{i + k} y^{i + j} z^{j + k}  \cdot
  $$ $$F\left( {\alpha  + i + j + k,\beta _1  + i + j + k;\gamma _1  + i + k;x}
  \right)F\left( {\beta _1  + i + j,\beta _2  + i + j;\gamma _2  + i + j;y} \right)$$
$$ F\left( {\beta _2  + j + k,\alpha  + i + j + k;\gamma _3  + j + k;z} \right),\eqno
(4.11)$$ $$ H_C \left( {\alpha ,\beta _1 ,\beta _2 ;\gamma
\,;x,y,z} \right) = $$ $$= \sum\limits_{i,j,k,l,r = 0}^\infty  {}
\frac{{\left( { - 1} \right)^{i + j} \left(
  \alpha  \right)_{i + 2j + k + l + r} \left( {\beta _1 } \right)_{2i + j + k + l}
  \left( {\beta _1 } \right)_{2i + j + k + r} \left( {\beta _2 } \right)_{i + j +
  k + l + r} }}{{\left( {\gamma  + i + j - 1} \right)_{i + j} \left( {\gamma  + 2i
  + 2j + k + l + r - 1} \right)_r \left( {\beta _1 } \right)_{2i + j + k} }}
  $$ $$ \frac{{\left( \gamma  \right)_{2i + 2j} \left( {\gamma  - \beta _2 }
  \right)_{i + j + r} }}{{\left( {\gamma  - \beta _2 } \right)_{i + j} \left(
   \gamma  \right)_{2i + 2j + k + l} \left( \gamma  \right)_{2i + 2j + k + l
   + 2r} i!j!k!l!r!}}x^{i + j + k + l} y^{i + k + r} z^{j + l + r}  \cdot $$
$$\cdot F\left( {\alpha  + i + 2j + k + l,\beta _1  + 2i + j + k + l;\gamma
+ 2i + 2j + k + l\,;x} \right) \cdot  $$ $$\cdot F\left( {\beta _2  + i + j + k +
 l + r,\beta _1  + 2i + j + k + r;\gamma  + 2i + 2j + k + l + 2r;y} \right)
 \cdot$$ $$\cdot F\left( {\beta _2  + i + j + k + l + r,\alpha  + i + 2j + k + l
 + r;\gamma  + 2i + 2j + k + l + 2r;z} \right) ,\eqno (4.12)$$ $$ H_A \left( {\alpha ,
 \beta _1 ,\beta _2 ;\gamma _1 ,\gamma _2 ;x,y,y} \right) = \sum\limits_{i,j =
 0}^\infty  {} \frac{{\left( \alpha  \right)_{i + j} \left( {\beta _1 } \right)_{i +
  j} \left( {\beta _2 } \right)_{i + j} }}{{\left( {\gamma _1 } \right)_{i + j}
  \left( {\gamma _2 } \right)_{i + j} i!j!}}x^{i + j} y^{i + j}  \cdot
  $$ $$\cdot F\left( {\alpha  + i + j,\beta _1  + i + j;\gamma _1  + i + j;x}
   \right)F\left( {\beta _2  + i + j;\alpha  + \beta _1  + 2i + j;\gamma _2
    + i + j;y} \right) \eqno (4.13)$$ $$H_B \left( {\alpha ,\beta _1 ,\beta _2 ;
    \gamma _1 ,\gamma ,\gamma ;x,y,y} \right) = \sum\limits_{i,j,k = 0}^\infty
     {} \frac{{\left( \alpha  \right)_{i + j + k} \left( {\beta _1 } \right)_{i
      + j + k} \left( {\beta _2 } \right)_{2i + j + k} }}{{\left( {\gamma _1 }
      \right)_{j + k} \left( \gamma  \right)_i \left( \gamma  \right)_{2i + j +
       k} i!j!k!}}x^{j + k} y^{2i + j + k}  $$ $$F\left( {\alpha  + i + j + k,
       \beta _1  + i + j + k;\gamma _1  + j + k;x} \right) $$ $$ F\left( {\beta _2
        + 2i + j + k;\alpha  + \beta _1  + 2i + 2j + k;\gamma  + 2i + j + k;y} \right)
        \eqno    (4.14)$$ $$H_C \left( {\alpha ,\beta _1 ,\beta _2 ;\gamma \,;x,y,y}
\right) =$$ $$ = \sum\limits_{i,j,k,l = 0}^\infty  {} \frac{{\left( { - 1}
\right)^{i + j} \left( \alpha  \right)_{i + 2j + k + l} \left( {\beta _1 }
\right)_{2i + j + k + l} \left( {\beta _2 } \right)_{i + j + k + l} \left(
 \gamma  \right)_{2i + 2j} }}{{\left( {\gamma  + i + j - 1} \right)_{i + j}
  \left( \gamma  \right)_{2i + 2j + k + l} \left( \gamma  \right)_{2i + 2j
   + k + l} i!j!k!l!}}x^{i + j + k + l} y^{i + j + k + l}  \cdot
   $$ $$\cdot F\left( {\alpha  + i + 2j + k + l,\beta _1  + 2i + j +
   k + l;\gamma  + 2i + 2j + k + l\,;x} \right)$$ $$F\left( {\beta _2  + i + j +
   k + l;\alpha  + \beta _1  + 3i + 3j + 2k + l;\gamma  + 2i + 2j + k + l\,;y}
   \right)\eqno    (4.15)$$\\[2mm]

{\bf 5. Proofs of the received decompositions.}

\medskip

{\bf Case 1.} Decomposition (4.1) can be proved with the help of
integral representation, i.e.:
$$H_A \left( {\alpha ,\beta _1 ,\beta _2 ;\gamma _1 ,\gamma _2 ;x,y,z} \right)
 = \frac{{\Gamma \left( {\gamma _1 } \right)\Gamma \left( {\gamma _2 } \right)}}
 {{\Gamma \left( {\beta _1 } \right)\Gamma \left( {\beta _2 } \right)\Gamma
 \left( {\gamma _1  - \beta _1 } \right)\Gamma \left( {\gamma _2  - \beta _2
 } \right)}} \cdot  $$ $$\cdot \int\limits_0^1 {\int\limits_0^1 {} }
 \xi ^{\beta _1  - 1} \eta ^{\beta _2  - 1} \left( {1 - \xi } \right)^{\gamma _1
  - \beta _1  - 1} \left( {1 - \eta } \right)^{\gamma _2  - \beta _2  - 1}
  \left( {1 - y\eta } \right)^{\alpha  - \beta _1 } \left[ {\left( {1 - y\eta }
   \right)\left( {1 - z\eta } \right) - x\xi } \right]^{ - \alpha } d\xi \,d\eta ,
\eqno  (5.1)$$ $${\mathop{\rm Re}\nolimits} \beta _1  >
0,\,\,{\mathop{\rm Re}\nolimits} \beta _2  > 0,\,\,{\mathop{\rm
Re}\nolimits} \left( {\gamma _1  - \beta _1 } \right) >
0,\,\,{\mathop{\rm Re}\nolimits} \left( {\gamma _2  - \beta _2 }
\right) > 0.$$

We shall prove decomposition (4.1). For this purpose, we shall
take advantage of integral representation (5.1). Taking into
account formula $$ \left[ {\left( {1 - y\eta } \right)\left( {1 -
z\eta } \right) - x\xi } \right]^{ - \alpha }  = \left[ {\left( {1
- x\xi } \right)\left( {1 - y\eta } \right)\left( {1 - z\eta }
\right)} \right]^{ - \alpha } \sum\limits_{i,j = 0}^\infty  {}
\frac{{\left( \alpha  \right)_{i + j} }}{{i!j!}}\sigma _1^i \sigma
_2^j,\eqno (5.2)$$ where $$ \sigma _1  = \frac{{xz\xi \eta
}}{{\left( {1 - x\xi } \right)\left( {1 - y\eta } \right)\left( {1
- z\eta } \right)}},\,\,\,\,\, \sigma _2  = \frac{{xy\xi \eta
}}{{\left( {1 - x\xi } \right)\left( {1 - y\eta } \right)}} .
\eqno (5.3)$$

Substituting equality (4.2) into integral representation (4.1), we
find
$$ H_A \left( {\alpha ,\beta _1 ,\beta _2 ;\gamma _1 ,\gamma _2 ;x,y,z} \right)
 = \sum\limits_{i,j = 0}^\infty  {} \frac{{\left( \alpha  \right)_{i + j} }}
 {{i!j!}}x^{i + j} y^j z^i $$ $$\cdot \frac{{\Gamma \left( {\gamma _1 } \right)}}
 {{\Gamma \left( {\beta _1 } \right)\Gamma \left( {\gamma _1  - \beta _1 }
 \right)}}\int\limits_0^1 {\,\xi ^{\beta _1  + i + j - 1} \left( {1 - \xi }
 \right)^{\gamma _1  - \beta _1  - 1} \left( {1 - x\xi } \right)^{ - \alpha
  - i - j} \,d\xi \,\,} $$ $$ \cdot \frac{{\Gamma \left( {\gamma _2 } \right)}}
  {{\Gamma \left( {\beta _2 } \right)\Gamma \left( {\gamma _2  - \beta _2 }
  \right)}}\int\limits_0^1 {} \eta ^{\beta _2  + i + j - 1} \left( {1 - \eta }
  \right)^{\gamma _2  - \beta _2  - 1} \left( {1 - y\eta } \right)^{ - \beta _1
   - i - j} \left( {1 - z\eta } \right)^{ - \alpha  - i} d\eta  .\eqno (5.4)$$

By virtue of integral representations $$\int\limits_0^1 {\,\xi ^{b
- 1} \left( {1 - \xi } \right)^{c - b - 1} \left( {1 - x\xi }
\right)^{ - a} \,d\xi \,\,}  = \frac{{\Gamma \left( b
\right)\Gamma \left( {c - b} \right)}}{{\Gamma \left( c
\right)}}F\left( {a,b;c;x} \right),$$ $${\mathop{\rm Re}\nolimits}
\,b > 0,\,{\mathop{\rm Re}\nolimits} \,\left( {c - b} \right) > 0
 ,$$ $$\int\limits_0^1 {} \eta ^{a - 1} \left( {1 - \eta } \right)^{c - a
- 1} \left( {1 - y\eta } \right)^{ - b_1 } \left( {1 - z\eta }
\right)^{ - b_2 } d\eta  = \frac{{\Gamma \left( a \right)\Gamma
\left( {c - a} \right)}}{{\Gamma \left( c \right)}}F_1 \left(
{a,b_1 ,b_2 ;c;y,z} \right),$$ $${\mathop{\rm Re}\nolimits} \,a
> 0,\,{\mathop{\rm Re}\nolimits} \,\left( {c - a} \right) > 0
,$$ from (4.4) follows decomposition (3.1).

{\bf Case 2.} The obtained operator identities (2.1) - (2.15)
consist of the operator or superposition of operators (1.1) -
(1.12). We shall show, as it is possible to apply superpositions
of operators for Hypergeometric function. For instance, we
consider decomposition (3.7). It's easy to see, that equality
takes place
$$ \nabla _{xz} \left( \alpha  \right)\nabla _{xy} \left( {\beta _1 }
\right)\nabla _{yz} \left( {\beta _2 } \right) = $$
  $$= \frac{1}{{\left( {\beta _2 } \right)_n \left( {\beta _2 } \right)_p
  \left( \alpha  \right)_p }}\sum\limits_{i,j,k = 0}^\infty  {} \frac{{
  \left( {\beta _2 } \right)_{n + j} \left( {\beta _2 } \right)_{p + i}
  \left( \alpha  \right)_{p + i} \left( { - \delta _1 } \right)_{i + j}
  \left( { - \delta _2 } \right)_{i + k} \left( { - \delta _3 } \right)_{j + k}
  }}{{\left( \alpha  \right)_{i + j} \left( {\beta _1 } \right)_i \left( {\beta _2
  } \right)_{i + j + k} i!j!k!}},\eqno (5.5)$$

Taking into account the identities (4.5), from parity (2.7), we
have
$$H_A \left( {\alpha ,\beta _1 ,\beta _2 ;\gamma _1 ,\gamma _2 ;x,y,z} \right)
= \sum\limits_{i,j,k = 0}^\infty  {} \frac{{\left( {\beta _2 } \right)_j \left(
{\beta _2 } \right)_i \left( \alpha  \right)_i \left( { - \delta _1 } \right)_{i
+ j} \left( { - \delta _2 } \right)_{i + k} \left( { - \delta _3 } \right)_{j + k}
}}{{\left( \alpha  \right)_{i + j} \left( {\beta _1 } \right)_i \left( {\beta _2 }
 \right)_{i + j + k} i!j!k!}} \cdot $$ $$\cdot F\left( {\alpha ,\beta _1 ;\gamma _1
 ;x} \right)F_3 \left( {\beta _1 ,\beta _2  + i,\beta _2  + j,\alpha  + i;\gamma _2
  ;y,z} \right) .\eqno  (5.6)$$

By virtue of the formula [16, p. 93] $$ \left( {\delta  + a}
\right)\left( {\delta  + a + 1} \right) \cdot \cdot  \cdot \left(
{\delta  + a + r - 1} \right)f\left( \xi \right) = \xi ^{1 - a}
\frac{{d^r }}{{d\xi ^{\,\,r} }}[\xi ^{a + r - 1} f\left( \xi
\right)],\eqno (5.7)$$

where $f\left( x \right)$ - analytic function, we find that
$$ \left( { - \delta } \right)_r f\left( \xi  \right) = \left( { -
1} \right)^r \xi ^r \frac{{d^r }}{{d\xi ^{\,\,r} }}f\left( \xi
\right).$$

We have
$$ \left( { - \delta _1 } \right)_{i + j} F\left( {\alpha ,\beta _1
;\gamma _1 ;x} \right) = \left( { - 1} \right)^{i + j}
\frac{{\left( \alpha  \right)_{i + j} \left( {\beta _1 }
\right)_{i + j} }}{{\left( {\gamma _1 } \right)_{i + j} }}x^{i +
j} F\left( {\alpha  + i + j,\beta _1  + i + j;\gamma _1  + i +
j;x} \right),\eqno (5.8)$$
$$ \left( { - \delta _2 } \right)_{i + k} \left( { - \delta _3 } \right)_{j
+ k} F_3 \left( {\beta _1 ,\beta _2  + i,\beta _2  + j,\alpha  + i;\gamma _2
 ;y,z} \right) = $$ $$= \left( { - 1} \right)^{i + j} y^{i + k} z^{j + k}
  \frac{{\left( \alpha  \right)_{i + j + k} \left( {\beta _1 } \right)_{i
   + k} \left( {\beta _2 } \right)_{i + j + k}^2 }}{{\left( \alpha  \right)_i
   \left( {\beta _1 } \right)_j \left( {\beta _2 } \right)_i \left( {\gamma _2
    } \right)_{i + j + 2k} }} \cdot  $$ $$\cdot F_3 \left( {\beta _1  + i + k,
    \beta _2  + i + j + k,\beta _2  + i + j + k,\alpha  + i + j + k;\gamma _2
     + i + j + 2k;y,z} \right) ,\eqno (5.9)$$ Substituting identities (4.8) -
(4.9) into equality (4.6), we obtain the decomposition (3.7).

{\bf Case 3. The integral representation exists for Hypergeometric
function $ H_B \left( {\alpha ,\beta _1 ,\beta _2 ;\gamma _1
,\gamma _2 ,\gamma _3 ;x,y,z} \right).$}

$$H_B \left( {\alpha ,\beta _1 ,\beta _2 ;\gamma _1 ,\gamma _2 ,
\gamma _3 ;x,y,z} \right) = \frac{{\Gamma \left( {\gamma _1 }
\right) \Gamma \left( {\gamma _2 } \right)\Gamma \left( {\gamma _3
} \right)}} {{\Gamma \left( {\beta _1 } \right)\Gamma \left(
{\beta _2 } \right) \Gamma \left( \alpha  \right)\Gamma \left(
{\gamma _1  - \beta _1 } \right)\Gamma \left( {\gamma _2  - \beta
_2 } \right)\Gamma \left( {\gamma _3  - \alpha } \right)}} \cdot
$$ $$\cdot \int\limits_0^1 {\int\limits_0^1 {\int\limits_0^1 {} } }
\xi ^{\beta _1  - 1} \eta ^{\beta _2  - 1} \zeta ^{\alpha  - 1}
\left( {1 - \eta } \right)^{\gamma _2  - \beta _2  - 1} \left( {1
- x\xi } \right)^{1 + \beta _2  - \gamma _3 } \left( {1 - y\eta }
\right)^{2 + \alpha  - \gamma _1  - \gamma _3 }\cdot$$ $$ \left(
{1 - x\xi  - z\zeta } \right)^{1 + \beta _1  - \beta _2  - \gamma
_1 }  \cdot\left[ {\left( {1 - \zeta } \right)\left( {1 - x\xi }
\right) \left( {1 - y\eta } \right) + xy\xi \eta \zeta }
\right]^{\gamma _3  - \alpha  - 1} \cdot$$ $$\cdot\left[ {\left(
{1 - \xi } \right)\left( {1 - y\eta } \right) \left( {1 - x\xi  -
z\zeta } \right) + yz\xi \eta \zeta } \right]^{\gamma _1
 - \beta _1  - 1} d\xi d\eta d\zeta  \eqno (5.10)$$
$$ {\mathop{\rm Re}\nolimits} \alpha  > 0,{\mathop{\rm
Re}\nolimits} \beta _1  > 0,\,\,{\mathop{\rm Re}\nolimits} \beta
_2  > 0,\,\,{\mathop{\rm Re}\nolimits} \left( {\gamma _1  - \beta
_1 } \right) > 0,\,\,{\mathop{\rm Re}\nolimits} \left( {\gamma _2
- \beta _2 } \right),\,\,{\mathop{\rm Re}\nolimits} \left( {\gamma
_3  - \alpha } \right) > 0 .$$

We shall prove the decomposition (3.11).
$$ \left[ {\left( {1 - \zeta } \right)\left( {1 - x\xi } \right)
\left( {1 - y\eta } \right) + xy\xi \eta \zeta } \right]^{\gamma _3
  - \alpha  - 1}  = \left( {1 - x\xi } \right)^{\gamma _3  -
  \alpha  - 1} \left( {1 - y\eta } \right)^{\gamma _3  - \alpha  - 1}
\left( {1 - \zeta } \right)^{\gamma _3  - \alpha  - 1}\cdot$$
$$\cdot \sum\limits_{i =
 0}^\infty  {} \frac{{\left( {1 + \alpha  - \gamma _3 } \right)_i }}{{i!}}
 \left( { - \frac{{xy\xi \eta \zeta }}{{\left( {1 - x\xi } \right)\left( {1
 - y\eta } \right)\left( {1 - \zeta } \right)}}} \right)^i  ,\eqno  (5.11)$$
$$ \left[ {\left( {1 - \xi } \right)\left( {1 - y\eta } \right)\left( {1 - x\xi
 - z\zeta } \right) + yz\xi \eta \zeta } \right]^{\gamma _1  - \beta _1  - 1}
 = \left( {1 - \xi } \right)^{\gamma _1  - \beta _1  - 1} \left( {1 - y\eta }
  \right)^{\gamma _1  - \beta _1  - 1}  \cdot $$ $$ \cdot \left( {1 - x\xi  -
  z\zeta } \right)^{\gamma _1  - \beta _1  - 1} \sum\limits_{j = 0}^\infty  {}
   \frac{{\left( {1 + \beta _1  - \gamma _1 } \right)}}{{j!}}\left( { -
   \frac{{yz\xi \eta \zeta }}{{\left( {1 - \xi } \right)\left( {1 - y\eta }
   \right)\left( {1 - x\xi  - z\zeta } \right)}}} \right)^j  .\eqno (5.12)$$

Substituting (4.11)-(4.12) in integral representation (4.10), we
define
$$ H_B \left( {\alpha ,\beta _1 ,\beta _2 ;\gamma _1 ,\gamma _2 ,\gamma _3
 ;x,y,z} \right) = \frac{{\Gamma \left( {\gamma _1 } \right)\Gamma \left(
 {\gamma _2 } \right)\Gamma \left( {\gamma _3 } \right)}}{{\Gamma \left(
 {\beta _1 } \right)\Gamma \left( {\beta _2 } \right)\Gamma \left( \alpha
  \right)\Gamma \left( {\gamma _1  - \beta _1 } \right)\Gamma \left( {\gamma _2
   - \beta _2 } \right)\Gamma \left( {\gamma _3  - \alpha } \right)}} \cdot $$
$$  \cdot \sum\limits_{i,j = 0}^\infty  {} \frac{{\left( { - 1} \right)^{i + j}
 \left( {1 + \alpha  - \gamma _3 } \right)_i \left( {1 + \beta _1  - \gamma _1 }
  \right)_j }}{{i!j!}}x^i y^{i + j} z^j  \cdot $$
$$  \cdot \int\limits_0^1 {\int\limits_0^1 {\int\limits_0^1 {} } } \xi ^{\beta _1
 - 1 + i + j} \eta ^{\beta _2  - 1 + i + j} \zeta ^{\alpha  - 1 + i + j}
 \left( {1 - \xi } \right)^{\gamma _1  - \beta _1  - 1 - j} \left( {1 -
  \eta } \right)^{\gamma _2  - \beta _2  - 1} \left( {1 - \zeta }
  \right)^{\gamma _3  - \alpha  - 1 - i}  \cdot  $$
$$  \cdot \left( {1 - x\xi } \right)^{\beta _2  - \alpha  - i} \left(
{1 - y\eta } \right)^{ - \beta _1  - i - j} \left( {1 - x\xi  -
z\zeta } \right)^{ - \beta _2  - j} d\xi d\eta d\zeta  .\eqno
(5.13)$$ By virtue of parity $$ \left( {1 - x\xi  - z\zeta }
\right)^{ - \beta _2  - j}  = \left( {1 - x\xi } \right)^{ - \beta
_2  - j} \left( {1 - z\zeta } \right)^{ - \beta _2  - j}
\sum\limits_{k = 0}^\infty  {} \frac{{\left( {\beta _2  + j}
\right)_k }}{{k!}}\left( {\frac{{xz\xi \zeta }}{{\left( {1 - x\xi
} \right)\left( {1 - z\zeta } \right)}}} \right)^k,\eqno (5.14)$$
we have $$ H_B \left( {\alpha ,\beta _1 ,\beta _2 ;\gamma _1
,\gamma _2 , \gamma _3 ;x,y,z} \right) $$
$$ = \sum\limits_{i,j,k = 0}^\infty  {} \frac{{\left( { - 1} \right)^{i
+ j} \left( {1 + \alpha  - \gamma _3 } \right)_i \left( {1 + \beta _1  -
 \gamma _1 } \right)_j \left( {\beta _2  + j} \right)_k }}{{i!j!k!}}x^{i
 + k} y^{i + j} z^{j + k}  $$
$$\cdot \frac{{\Gamma \left( {\gamma _1 } \right)}}{{\Gamma \left( {\beta _1 }
 \right)\Gamma \left( {\gamma _1  - \beta _1 } \right)}}\int\limits_0^1 {}
 \xi ^{\beta _1  - 1 + i + j + k} \left( {1 - \xi } \right)^{\gamma _1  -
 \beta _1  - 1 - j} \left( {1 - x\xi } \right)^{ - \alpha  - i - j - k} d\xi $$
$$  \cdot \frac{{\Gamma \left( {\gamma _2 } \right)}}{{\Gamma \left( {\beta _2 }
\right)\Gamma \left( {\gamma _2  - \beta _2 }
\right)}}\int\limits_0^1 {} \eta ^{\beta _2  - 1 + i + j} \left(
{1 - \eta } \right)^{\gamma _2  - \beta _2  - 1} \left( {1 - y\eta
} \right)^{ - \beta _1  - i - j} d\eta $$ $$
  \cdot \frac{{\Gamma \left( {\gamma _3 } \right)}}{{\Gamma \left( \alpha
   \right)\Gamma \left( {\gamma _3  - \alpha } \right)}}\int\limits_0^1 {}
    \zeta ^{\alpha  - 1 + i + j + k} \left( {1 - \zeta } \right)^{\gamma _3
     - \alpha  - 1 - i} \left( {1 - z\zeta } \right)^{ - \beta _2  - j - k}
     d\zeta .  \eqno (5.15)$$ and from identity (5.15) we define the decomposition (3.11).

{\bf Case 4.} Decomposition (3.8) can be proved by means of
equality
$$ \nabla _{xz} \left( \alpha  \right)\nabla _{xy} \left( {\beta _1 } \right)
\nabla _{yz} \left( {\beta _2 } \right)\nabla _{yz} \left( \gamma  \right) =
\frac{1}{{\left( \alpha  \right)_m \left( \alpha  \right)_p \left( {\beta _1 }
 \right)_m \left( {\beta _2 } \right)_n \left( {\beta _2 } \right)_p }} $$
$$  \cdot \sum\limits_{i,j,k,l = 0}^\infty  {\frac{{\left( \alpha  \right)_{j
+ k} \left( \alpha  \right)_k \left( {\beta _1 } \right)_{i + j} \left( {\beta _2 }
 \right)_i^2 \left( {\alpha  + j + k} \right)_m \left( {\alpha  + k} \right)_p
 \left( {\beta _1  + i + j} \right)_m \left( {\beta _2  + i} \right)_n
 \left( {\beta _2  + i} \right)_p }}{{\left( \alpha  \right)_{i + j + k +
 l} \left( {\beta _1 } \right)_{i + j + k} \left( {\beta _2 } \right)_{2i
 + j} \left( \gamma  \right)_i }}} $$ $$ \cdot \frac{{\left( { - \delta _1 }
  \right)_{k + l} \left( { - \delta _2 } \right)_{i + j + k} \left( { -
  \delta _3 } \right)_{i + j + l} }}{{i!j!k!l!}} \eqno  (5.16)$$

Taking into account the identities (4.16), from parity (2.8), we
have
$$ H_B \left( {\alpha ,\beta _1 ,\beta _2 ;\gamma _1 ,\gamma ,\gamma ;x,y,z}
 \right) = \sum\limits_{i,j,k,l = 0}^\infty  {\frac{{\left( \alpha  \right)_{j
 + k} \left( \alpha  \right)_k \left( {\beta _1 } \right)_{i + j}
 \left( {\beta _2 } \right)_i^2 }}{{\left( \alpha  \right)_{i + j +
 k + l} \left( {\beta _1 } \right)_{i + j + k} \left( {\beta _2 }
 \right)_{2i + j} \left( \gamma  \right)_i i!j!k!l!}}}  \cdot $$
$$  \cdot \left( { - \delta _1 } \right)_{k + l} F\left( {\alpha  + j
+ k,\beta _1  + i + j;\gamma _1 ;x} \right) \cdot  $$
$$  \cdot \left( { - \delta _2 } \right)_{i + j + k} \left( { -
\delta _3 } \right)_{i + j + l} F_3 \left( {\beta _2  + i,\alpha
 + k,\beta _1 ,\beta _2  + i;\gamma ;y,z} \right) \eqno (5.17)$$

We have
$$ \left( { - \delta _1 } \right)_{k + l} F\left( {\alpha  + j + k,\beta _1
+ i + j;\gamma _1 ;x} \right) =  $$ $$ = \left( { - 1} \right)^{k + l} x^{k +
l} \frac{{\left( \alpha  \right)_{j + 2k + l} \left( {\beta _1 } \right)_{i +
j + k + l} }}{{\left( \alpha  \right)_{j + k} \left( {\beta _1 } \right)_{i +
j} \left( {\gamma _1 } \right)_{k + l} }}F\left( {\alpha  + j + 2k + l,\beta _1
  + i + j + k + l;\gamma _1  + k + l;x} \right) ,\eqno (5.18)$$
$$ \left( { - \delta _2 } \right)_{i + j + k} \left( { - \delta _3 } \right)_{i
+ j + l} F_3 \left( {\beta _2  + i,\alpha  + k,\beta _1 ,\beta _2
+ i;\gamma ; y,z} \right) =  $$ $$
  = \left( { - 1} \right)^{k + l} y^{i + j + k} z^{i + j + l} \frac{{\left(
   \alpha  \right)_{i + j + k + l} \left( {\beta _1 } \right)_{i + j + k}
   \left( {\beta _2 } \right)_{2i + j + k} \left( {\beta _2 } \right)_{2i +
    j + l} }}{{\left( \alpha  \right)_k \left( {\beta _2 } \right)_i^2
    \left( \gamma  \right)_{2i + 2j + k + l} }} $$
$$ F_3 ( \beta _2  + 2i + j + k,\alpha  + i + j + k + l,\beta _1
 + i + j + k,\beta _2  + 2i + j + l;$$ $$\gamma  + 2i + 2j + k + l;y,z )\eqno  (5.19)$$

Substituting identities (4.18) - (4.19) into equality (4.17), we
obtain the decomposition (3.8).\\[2mm]

{\bf 6. The Integrals.}
\medskip

With the help of the obtained decompositions can be easily found
the following integrals connecting with Hypergeometric  $H_A
,\,\,H_B $ functions
$$ H_A \left( {\alpha ,\beta _1 ,\beta _2 ;\gamma _1 ,\gamma _2 ;x,y,z}
\right) = \frac{{\Gamma \left( {\gamma _1 } \right)\Gamma \left( {\gamma _2
} \right)}}{{\Gamma \left( {\beta _1 } \right)\Gamma \left( {\beta _2 } \right)
\Gamma \left( {\gamma _1  - \beta _1 } \right)\Gamma \left( {\gamma _2  -
 \beta _2 } \right)}} \cdot  $$
$$  \cdot \int\limits_0^1 {\int\limits_0^1 {\xi ^{\beta _1  - 1} \eta ^{\beta _2
  - 1} \left( {1 - \xi } \right)^{\gamma _1  - \beta _1  - 1} \left( {1 - \eta }
   \right)^{\gamma _2  - \beta _2  - 1} \left( {1 - y\eta } \right)^{ - \beta _2 }
   \left( {1 - x\xi  - z\eta } \right)^{ - \alpha }  \cdot } } $$
$$ \cdot F\left( {\alpha ,\beta _2 ;\beta _1 ;\frac{{xy\xi \eta }}{{\left(
{1 - y\eta } \right)\left( {1 - x\xi  - z\eta } \right)}}}
\right)d\xi \,d\eta  ,\eqno (6.1)$$ $${\mathop{\rm Re}\nolimits}
\alpha _1  > 0,\,\,\,\,{\mathop{\rm Re}\nolimits} \left( {\gamma
_1  - \alpha _1 } \right) > 0,\,\,\,\,\,\,\,{\mathop{\rm
Re}\nolimits} \alpha _2  > 0,\,\,\,\,{\mathop{\rm Re}\nolimits}
\left( {\gamma _2  - \alpha _2 } \right) > 0,$$
$$ H_B \left( {\alpha ,\beta _1 ,\beta _2 ;\gamma _1 ,\gamma _2 ,\gamma _3
;x,y,z} \right) = \frac{{\Gamma \left( {\gamma _1 } \right)\Gamma \left(
{\gamma _2 } \right)\Gamma \left( {\gamma _3 } \right)}}{{\Gamma ^2
\left( \alpha  \right)\Gamma \left( {\beta _1 } \right)\Gamma \left(
 {\gamma _1  - \alpha } \right)\Gamma \left( {\gamma _2  - \beta _1 }
 \right)\Gamma \left( {\gamma _3  - \alpha } \right)}} \cdot  $$
$$  \cdot \int\limits_0^1 {\int\limits_0^1 {\int\limits_0^1 {} } } \xi ^{\alpha
 - 1} \eta ^{\beta _1  - 1} \zeta ^{\alpha  - 1} \left( {1 - \xi }
 \right)^{\gamma _1  - \alpha  - 1} \left( {1 - \eta } \right)^{\gamma _2
  - \beta _1  - 1} \left( {1 - \zeta } \right)^{\gamma _3  - \alpha  - 1}
  \left( {1 - x\xi } \right)^{\beta _2  - \beta _1 }  \cdot  $$
$$  \cdot \left[ {\left( {1 - x\xi } \right)\left( {1 - y\eta  - z\zeta }
 \right) - xy\xi \eta } \right]^{ - \beta _2 }  \cdot  $$
$$  \cdot F\left( {\beta _2 ,1 + \beta _1  - \gamma _2 ;\alpha ; -
\frac{{xz\xi \eta \zeta }}{{\left( {1 - \eta } \right)\left[ {\left(
{1 - x\xi } \right)\left( {1 - y\eta  - z\zeta } \right) - xy\xi \eta }
 \right]}}} \right)d\xi d\eta d\zeta  ,\eqno   (6.2)$$
$${\mathop{\rm Re}\nolimits} \alpha  > 0,\,\,\,\,{\mathop{\rm
Re}\nolimits} \left( {\gamma _1  - \alpha } \right) >
0,\,{\mathop{\rm Re}\nolimits} \beta _1  > 0,\,\,\,\,{\mathop{\rm
Re}\nolimits} \left( {\gamma _2  - \beta _1 } \right) >
0,\,\,\,\,{\mathop{\rm Re}\nolimits} \left( {\gamma _3  - \alpha }
\right) > 0.$$

These identities can be proved by using the expansions of under
integral functions.

If we use the definition for Euler function $\Gamma \left( z
\right)$ ([13], Ch. 1, Sec. 1.2, (1))

$$\Gamma \left( z \right) = \left\{ {\begin{array}{*{20}c}
   {\int\limits_0^\infty  {} e^{ - t} t^{z - 1} dt,\,{\mathop{\rm Re}\nolimits}
   \,\left( z \right) > 0,}  \\
   {\frac{{\Gamma \left( {1 + z} \right)}}{z},{\mathop{\rm Re}\nolimits}
   \,\left( z \right) < 0;\,z \ne  - 1, - 2,....}  \\
\end{array}} \right.\eqno(6.3)$$

Then, for Hypergeometric $H_A ,\,H_B ,\,H_C $ functions we find
the following integral representations

$$ H_A \left( {\alpha ,\beta _1 ,\beta _2 ;\gamma _1 ,\gamma _2 ;x,y,z}
\right) = \frac{1}{{\Gamma \left( \alpha  \right)\Gamma \left(
{\beta _1 } \right)\Gamma \left( {\beta _2 } \right)}} $$
$$  \cdot \int\limits_0^\infty  {\int\limits_0^\infty  {\int\limits_0^\infty
{} } } e^{ - u_1  - u_2  - u_3 } u_1 ^{\alpha  - 1} u_2 ^{\beta _1
- 1} u_3 ^{\beta _2  - 1} {}_0F_1 \left( {\gamma _1 ;xu_1 u_2 }
\right){}_0F_1 \left( {\gamma _2 ;yu_2 u_3  + zu_1 u_3 }
\right)du_1 du_2 du_3 ,\eqno  (6.4)$$
$${\mathop{\rm Re}\nolimits} \,\alpha  > 0,\,{\mathop{\rm
Re}\nolimits} \,\beta _1  > 0,\,{\mathop{\rm Re}\nolimits} \,\beta
_2  > 0,$$
$$ H_A \left( {\alpha ,\beta _1 ,\beta _2 ;\gamma _1 ,\gamma _2 ;x,y,z} \right) $$
$$  = \frac{1}{{\Gamma \left( \alpha  \right)\Gamma \left( {\beta _1 } \right)}}
\int\limits_0^\infty  {\int\limits_0^\infty  {} } e^{ - u_1  - u_2 } u_1 ^{\alpha
 - 1} u_2 ^{\beta _1  - 1} {}_0F_1 \left( {\gamma _1 ;xu_1 u_2 } \right)\,{}_1F_1
 \left( {\beta _2 ;\gamma _2 ;yu_2  + zu_1 } \right)du_1 du_2 ,\eqno  (6.5)$$
$${\mathop{\rm Re}\nolimits} \,\alpha  > 0,\,{\mathop{\rm
Re}\nolimits} \,\beta _1  > 0,$$
$$ H_B \left( {\alpha ,\beta _1 ,\beta _2 ;\gamma _1 ,\gamma _2 ,\gamma _3
;x,y,z} \right) = \frac{1}{{\Gamma \left( \alpha  \right)\Gamma \left(
 {\beta _1 } \right)\Gamma \left( {\beta _2 } \right)}} $$
$$  \cdot \int\limits_0^\infty  {\int\limits_0^\infty  {\int\limits_0^\infty
 {} } } e^{ - \xi  - \eta  - \zeta } \xi ^{\alpha  - 1} \eta ^{\beta _1  -
  1} \zeta ^{\beta _2  - 1} {}_0F_1 \left( {\gamma _1 ;x\xi \eta } \right)
  \,{}_0F_1 \left( {\gamma _2 ;y\eta \zeta } \right)\,{}_0F_1 \left(
   {\gamma _3 ;z\xi \zeta } \right)d\xi d\eta d\zeta ,\eqno  (6.6)$$
$${\mathop{\rm Re}\nolimits} \,\alpha  > 0,\,{\mathop{\rm Re}\nolimits}
\,\beta _1  > 0,\,{\mathop{\rm Re}\nolimits} \,\beta _2  > 0 ,$$
$$ H_B \left( {\alpha ,\beta _1 ,\beta _2 ;\gamma _1 ,\gamma _2 ,\gamma _3
;x,y,z} \right) $$
$$= \frac{1}{{\Gamma \left( \alpha  \right)\Gamma \left( {\beta _1 } \right)}}
\int\limits_0^\infty  {\int\limits_0^\infty  {} } e^{ - \xi  -
\eta } \xi ^{\alpha  - 1} \eta ^{\beta _1  - 1} {}_0F_1 \left(
{\gamma _1 ;x\xi \eta } \right)\Psi _2 \left( {\beta _2 ;\gamma _2
,\gamma _3 ;y\eta ,z\xi } \right)d\xi d\eta  ,\eqno (6.7)$$
$${\mathop{\rm Re}\nolimits} \,\alpha  > 0,\,{\mathop{\rm
Re}\nolimits} \,\beta _1  > 0,$$
$$ H_C \left( {\alpha ,\beta _1 ,\beta _2 ;\gamma \,;x,y,z} \right) $$
$$  = \frac{1}{{\Gamma \left( \alpha  \right)\Gamma \left( {\beta _1 } \right)
\Gamma \left( {\beta _2 } \right)}}\int\limits_0^\infty  {\int\limits_0^\infty
 {\int\limits_0^\infty  {} } } e^{ - \xi  - \eta  - \zeta } \xi ^{\alpha  - 1}
 \eta ^{\beta _1  - 1} \zeta ^{\beta _2  - 1} {}_0F_1 \left( {\gamma ;x\xi \eta
  + y\eta \zeta  + z\xi \zeta } \right)d\xi d\eta d\zeta ,\eqno  (6.8)$$
  $${\mathop{\rm Re}\nolimits} \,\alpha  > 0,\,{\mathop{\rm Re}\nolimits}
  \,\beta _1  > 0,\,{\mathop{\rm Re}\nolimits} \,\beta _2  > 0,$$
where [15] $$\Psi _2 \left( {\alpha ;\gamma _1 ,\gamma _2 ;x,y}
\right) = \sum\limits_{m,n = 0}^\infty  {\frac{{\left( \alpha
\right)_{m + n} }}{{\left( {\gamma _1 } \right)_m \left( {\gamma
_2 } \right)_n m!n!}}x^m y^n },\eqno (6.9)$$ and, as usual $
{}_pF_q $ denotes the generalized Hypergeometric function with $
p$ numerator and $ q$  denominator parameters. For multiple
hypergeometric of functions $ H_B^{(n)} ,\,\,H_C^{(n)} $ the
integrated representations are determined in the monograph ([13],
Ch. 9, Sec. 9.4, (198)-(200)).\\[2mm]

{\bf 7. Conclusion.}

\medskip

From decomposions (4.1)-(4.15) as a private value can be found the
known decompositions that have been found in [5], i.e.
decompositions (4.1)-(4.15) generalize the results of [5].
Further, the following decompositions for Hypergeometric Appell's
functions
$$
 F_1 \left( {\alpha ;\beta ,\beta ;\gamma ;x,y} \right) = \sum\limits_{i,j
= 0}^\infty  {\frac{{\left( { - 1} \right)^{i + j} \left( \alpha
\right)_{2i + 2j} \left( \beta  \right)_i \left( \beta  \right)_{i
+ j} \left( \gamma \right)_{2i} }}{{\left( {\gamma  + i - 1}
\right)_i \left( \gamma  \right)_{2i + j}^2 i!j!}}} x^{i + j} y^{i
+ j}
$$
$$
F_4 \left( {\alpha  + 2i + 2j,\beta  + i + j;\gamma  +
2i + j,\gamma + 2i + j;x,y} \right) ,\eqno (7.1)$$
$$ F_2 \left( {\alpha ;\beta _1 ,\beta _2 ;\gamma ,\gamma ;x,y}
 \right) = \sum\limits_{i,j = 0}^\infty  {} \frac{{\left( \alpha
 \right)_{2i + j} \left( {\beta _1 } \right)_{i + j} \left( {\beta _2 }
  \right)_{i + j} }}{{\left( \gamma  \right)_i \left( \gamma  \right)_{2i
  + 2j} i!j!}}x^{i + j} y^{i + j}
$$
$$
F_3 \left( {\alpha  + 2i + j,\alpha  + 2i + j;\beta _1  + i + j,\beta _2
 + i + j;\gamma  + 2i + 2j;x,y} \right) ,\eqno  (7.2)
$$
have not been found in work [5]. In particular, from
decompositions (4.1), (4.2) we can find the following identities
$$
 H_A \left( {\alpha ,\beta _1
,\beta _2 ;\gamma _1 ,\beta _2 ;x,y,z} \right) = \left( {1 - y}
\right)^{ - \beta _1 } \left( {1 - z} \right)^{ - \alpha } F\left(
{\alpha ,\beta _1 ;\gamma _1 ;\frac{x}{{\left( {1 - y}
\right)\left( {1 - z} \right)}}} \right) ,\eqno (7.3)
$$
$$
H_B \left( {\alpha ,\beta _1 ,\beta _2 ;\gamma _1 ,\beta _2 ,\beta _2
;x,y,z} \right)
$$
$$
  = \left( {1 - y} \right)^{ - \beta _1 }
\left( {1 - z} \right)^{ - \alpha } F_4 \left( {\alpha ,\beta _1
;\beta _2 ,\gamma _1 ; \frac{{yz}}{{\left( {1 - y} \right)\left(
{1 - z} \right)}},\frac{x}{{\left( {1 - y} \right)\left( {1 - z}
\right)}}} \right).\eqno (7.4)
$$

It should be noted that the following superpositions of the
operators $$ \nabla _{xy} \left( \alpha  \right)\nabla _{xz}
\left( \alpha  \right)\Delta _{yz} \left( \alpha  \right)\Delta
_{yz} \left( {\gamma _2 } \right),\,\,\,\,\,\,\nabla _{xz} \left(
\alpha \right)\nabla _{xy} \left( \alpha \right)\Delta _{yz}
\left( \alpha  \right),$$ $$\nabla _{xz} \left( \alpha
\right)\nabla _{xy} \left( \alpha \right)\Delta _{yz} \left(
\alpha \right)\Delta _{yz} \left( \gamma  \right)\tilde \Delta _x
\left( \gamma  \right) $$ have not been applied for Hypergeometric
functions.

\begin{center}
 REFERENCES
\end{center}

[1] G.Lohofer. Theory of an electromagnetically levitated metal
spheres. Absorbed Power. SIAM J. Appl. Math., 1989, 49, N2,
567-581.

[2] A.W.Ninkkanen. Generalized hypergeometric series $ {}^N
F\left( {x_1 ,x_2 ,...,x_N } \right)$ , arising in physical and
quantum chemical applications. J. Phys. A: Math. and Gen., 1983,
16, N9, 1813-1825.

[3] F. I. Frankl. The elected works on gas dynamics. Moscow, 1973.

[4] A. Hasanov. About one mixed task for the equation $signy\left|
y \right|^m u_{xx}  + x^n u_{yy}  = 0.$ New. Academ. Scien. Uz
SSR, N 2, 1982, 28-32.

[5] J. L. Burchnall and T.W. Chaundy. Expansions of Appell's
double hypergeometric functions. Quart. J. Math., Oxford, 1940,
Ser. 11, 249-270.

[6] J. L. Burchnall and T.W. Chaundy. Expansions of Appell's
double hypergeometric functions (II). Quart. J. Math., Oxford,
1941, Ser. 12, 112-128.

[7] T.W. Chaundy. Expansions of hypergeometric functions. Quart.
J. Math., Oxford, 1942, Ser. 13, 159-171.

[8] Srivastava H.M., Hypergeometric functions of three variables.
Ganita, 1964, 15, N2, 97-108.

[9] J. P. Singhal and S. S. Bhati, Certain expansions associated
with hypergeometric functions of $ n$ variables, Glasnik Mat. Ser.
3, 11 (31), 1976, 239-245.

[10] H. M. Srivastava, Some integrals representing triple
hypergeometric functions, Rend. Circ. Mat. Palermo, (ser. 2) 16,
1967, 99-115.

[11] G. Lauricella, Sulle funzioni ipergeometriche a piu
variabili, Rend. Circ. Mat. Palermo 7, (1893), 111-158.

[12] P.
Appell and J. Kampe de Feriet, Fonctions hypergeometriques et
hyperspheriques. Polynomes d'Hermite, Gauthier - Villars. Paris,
1926.

[13] H. M. Srivastava and P. W. Karlsson, Multiple Gaussian
Hypergeometric Series. Halsted Press (Ellis Horwood Limited,
Chichester), John Wiley and Sons, New York, Chichester, Brisbane
and Toronto, 1985.

[14] O. I. Marichev. Handbook of integral transforms of higher
transcendental functions, theory and algorithmic tables.
Chichester, Ellis Horwood, 1982.

[15] A. Erdelyi (Ed): Higher transcendental functions, vol. 1. New
York: McGraw Hill Book Company, Inc.  1953.

[16] E. G. C. Poole. Introduction to the theory of linear
differential equation. Oxford, 1936.

\end{document}